\documentclass[aap,MSNbibl,dvips]{arximspdf}
\usepackage{graphicx}
\usepackage{url,breakurl}

\doi{10.1214/13-AAP993} %kopijuoti is PTS
\volume{25}
\issue{1}
\pubyear{2015}
\firstpage{211}
\lastpage{234}

\makeatletter
\newcommand{\Rsta}{\Rightarrow^{\stably}}%{\mathop{\Rightarrow}^{\stably}}
\newcommand{\MC}{\mathrm{MC}}
\newcommand{\MMC}{\mathrm{MMC}}
\newcommand{\stably}{\mathrm{stably}}
\newcommand{\Var}{\operatorname{Var}}
\newcommand{\rrvert}{\vert}
\newcommand{\llvert}{\vert}
\newtheorem{theorem}{Theorem}
\newtheorem{proposition}{Proposition}
\newtheorem{lemma}{Lemma}
\newproclaim{remark}{Remark}
\makeatother

\begin{document}
\begin{frontmatter}

\title{Central limit theorem for the multilevel Monte~Carlo Euler method}
\runtitle{The multilevel Monte Carlo method}

\begin{aug}
\author[A]{\fnms{Mohamed} \snm{Ben Alaya}\thanksref{T1}\ead[label=e1]{mba@math.univ-paris13.fr}\ead[label=u1,url]{http://labex-mme-dii.u-cergy.fr/}}
\and
\author[A]{\fnms{Ahmed} \snm{Kebaier}\corref{}\thanksref{T2}\ead[label=e2]{kebaier@math.univ-paris13.fr}\ead[label=u2,url]{http://labex-mme-dii.u-cergy.fr/}}
\runauthor{M. Ben Alaya and A. Kebaier}
\affiliation{Universit\'{e} Paris 13}
\address[A]{LAGA, CNRS (UMR 7539)\\
Universit\'{e} Paris 13\\
Sorbonne Paris Cit\'{e}\\
99, av. J.B. Cl\'ement 93430 Villetaneuse\\
France\\
\printead{e1}\\
\phantom{E-mail:\ }\printead*{e2}} %adresu isvedimo komanda gale!
\end{aug}
\thankstext{T1}{Supported by Laboratory of Excellence MME-DII.\newline \printead{u1}.}
\thankstext{T2}{This research benefited from the support of the chair
``Risques Financiers,''
Fondation du Risque. Also supported by Laboratory of Excellence
MME-DII. \printead{u2}.}

% HISTORY:
\received{\smonth{3} \syear{2013}}
\revised{\smonth{9} \syear{2013}}

% ABSTRACT
%
\begin{abstract}
This paper focuses on studying the multilevel Monte Carlo method
recently introduced by Giles [\textit{Oper. Res.} \textbf{56} (2008) 607--617]
which is
significantly more efficient than the classical Monte Carlo one.
Our aim is to prove a central limit theorem of Lindeberg--Feller type
for the multilevel Monte Carlo method associated with the Euler discretization
scheme. To do so,
we prove first a stable law convergence theorem, in the spirit of Jacod
and Protter [\textit{Ann. Probab.} \textbf{26} (1998) 267--307],
for the Euler scheme error on two consecutive levels of the algorithm.
This leads to an accurate description of the optimal choice of
parameters and to
an explicit characterization of the limiting variance in the central
limit theorem of the algorithm. A complexity of the multilevel Monte
Carlo algorithm is carried out.
\end{abstract}

% KEYWORDS
% Pirmas kwd is didziosios raides
%
\begin{keyword}[class=AMS]
\kwd{60F05}
\kwd{62F12}
\kwd{65C05}
\kwd{60H35}.
\end{keyword}
\begin{keyword}
\kwd{Central limit theorem}
\kwd{multilevel Monte Carlo methods}
\kwd{Euler scheme}
\kwd{finance}
\end{keyword}

\end{frontmatter}

%s1 #&#
\section{Introduction}\label{intro}
In many applications, in particular in the pricing of financial securities,
we are interested in the effective computation by Monte Carlo methods
of the quantity $\mathbb Ef(X_T)$, where $X:=(X_t)_{0\leq t\leq T}$ is a
diffusion process and $f$ a~given function. The Monte Carlo Euler
method consists of two steps. First, approximate the diffusion process
$(X_t)_{0\leq t\leq T}$ by the Euler scheme $(X^n_t)_{0\leq t\leq T}$
with time step $T/n$. Then\vspace*{1pt} approximate $ \mathbb E f (
X^n_T )$ by
$
\frac{1}{N}\sum_{i=1}^{N}f( X^{n}_{T,i})$,
where $f( X^{n}_{T,i})_{1\leq i\leq N}$ is a sample of $N$ independent
copies of $f(X^n_T)$. This approximation is affected, respectively,
by a discretization error and a statistical error
\[
\varepsilon_n:=\mathbb E \bigl( f\bigl(X^n_T
\bigr)-f(X_T) \bigr)\quad\mbox{and}\quad \frac{1}{N}\sum
_{i=1}^{N}f\bigl( X^{n}_{T,i}
\bigr)-\mathbb Ef\bigl(X^n_T\bigr).
\]
On one hand, Talay and Tubaro \cite{TT} prove that if $f$ is
sufficiently smooth, then $\varepsilon_n\sim c/n$ with $c$ a given
constant and
in a more general context, Kebaier \cite{Keb} proves that the rate of
convergence of the discretization error $\varepsilon_n$ can be
$1/n^{\alpha}$ for all
values of $\alpha\in[1/2,1]$ (see, e.g., Kloeden and Platen \cite{kloeden}
for more details on discretization schemes). On the other hand, the
statistical error is controlled by the central limit theorem
with order $1/\sqrt{N}$. Further, the optimal choice of the sample size
$N$ in the classical Monte Carlo method
mainly depends on the order of the discretization error. More
precisely, it turns out that for $\varepsilon_n=1/n^{\alpha}$ the
optimal choice of $N$ is $n^{2\alpha}$.
This leads to a total complexity in the Monte Carlo method of order
$C_{\MC}=n^{2\alpha+1}$ (see Duffie and Glynn \cite{IIII} for related results).
Let us recall that the complexity of an algorithm is proportional to the
maximum number of basic computations performed by this one. Hence,
expressing this complexity in terms of the discretization error
$\varepsilon_n$, we get $C_{\MC}=\varepsilon_n^{-2-1/\alpha}$.

In order to improve the performance of this method, Kebaier introduced
a two-level Monte Carlo method \cite{Keb} (called the statistical
Romberg method) reducing the complexity $C_{\MC}$ while maintaining the
convergence of the algorithm. This method uses two Euler schemes with
time steps $T/n$ and $T/n^{\beta}$, $\beta\in(0,1)$ and approximates
$\mathbb E f(X_T)$ by
\[
\frac{1}{N_1}\sum_{i=1}^{N_1} f\bigl(
\widehat X^{n^{\beta}}_{T,i}\bigr) + \frac{1}{N_2} \sum
_{i=1}^{N_2}f\bigl( X^{n}_{T,i}
\bigr)-f\bigl( X^{n^{\beta}}_{T,i}\bigr),
\]
where\vspace*{-1pt} $\widehat X^{n^\beta}_T$ is a second Euler scheme with time step
$T/n^{\beta}$ and such that the Brownian paths used for $X^n_T$
and $ X^{n^{\beta}}_T$ has to be independent of the Brownian paths
used to simulate $\widehat X^{n^{\beta}}_T$.
It turns out that for a given discretization error $\varepsilon
_n=1/n^\alpha$ ($\alpha\in[1/2,1]$), the optimal choice is obtained
for $\beta=1/2$, $N_1=n^{2\alpha}$ and $N_2=n^{2\alpha-(1/2)}$. With
this choice, the complexity of the statistical Romberg method is of
order $C_{\mathrm{SR}}=n^{2\alpha+(1/2)}=\varepsilon_n^{-2-1/(2\alpha)}$, which
is lower than the classical complexity in the Monte Carlo method.

More recently, Giles \cite{Gil} generalized the statistical Romberg
method of Kebaier \cite{Keb}
and proposed the multilevel Monte Carlo algorithm, in a similar
approach to Heinrich's multilevel method for parametric
integration \cite{hei} (see also Creutzig et~al. \cite{Cre}, Dereich
\cite{Der}, Giles \cite{GilMil},
Giles, Higham and Mao \cite{GilHig}, Giles and Szpruch \cite{GilSzp},
Heinrich \cite{Heibis}, Heinrich and Sindambiwe \cite{HeiSin} and
Hutzenthaler, Jentzen and Kloeden
\cite{HutJenKlo} for related results).
The multilevel Monte Carlo method uses information from a sequence of
computations with decreasing step sizes and
approximates the quantity $\mathbb Ef( X_T)$ by
\[
Q_n=\frac{1}{N_0}\sum_{k=1}^{N_0}f
\bigl( X^{1}_{T,k}\bigr)+\sum_{\ell=1}^{L}
\frac{1}{N_\ell} \sum_{k=1}^{N_\ell} \bigl(f
\bigl( X^{\ell,m^{\ell}}_{T,k}\bigr)-f\bigl( X^{\ell,m^{\ell
-1}}_{T,k}
\bigr) \bigr),\qquad
m\in\mathbb N\setminus\{0,1\},
\]
where\vspace*{-2pt} the fine discretization step is equal to $T/n$ thereby $L=\frac
{\log n}{\log m}$. For $\ell\in\{1,\ldots,L\}$,
processes $(X^{\ell,m^{\ell}}_{t,k},X^{\ell,m^{\ell
-1}}_{t,k})_{0\leq
t\leq T}$, $k\in\{1,\ldots,N_\ell\}$, are independent copies
of $(X^{\ell,m^{\ell}}_t,X^{\ell,m^{\ell-1}}_t)_{0\leq t\leq T}$ whose\vspace*{1pt}
components denote the Euler schemes with time steps $m^{-\ell}T$ and
$m^{-(\ell-1)}T$.
However, for fixed $\ell$, the simulation of $(X^{\ell,m^{\ell
}}_t)_{0\leq t\leq T}$ and $ (X^{\ell,m^{\ell-1}}_t)_{0\leq t\leq T}$
has to be based on the same Brownian path. Concerning the first
empirical mean, processes
$( X^{1}_{t,k})_{0\leq t\leq T}$, $k\in\{1,\ldots,N_0\}$, are
independent copies of $( X^{1}_{t})_{0\leq t\leq T}$ which denotes the
Euler scheme
with time step~$T$. Here, it is important to point out that all these
$L+1$ Monte
Carlo estimators have to be based on different independent samples.
Due to the above independence assumption on the paths, the
variance of the multilevel estimator is given by
\[
\sigma^2:=\Var(Q_n)=N^{-1}_{0}\Var
\bigl(f\bigl(X^{1}_{T}\bigr)\bigr)+\sum
_{\ell=1}^L N^{-1}_{\ell}
\sigma^2_{\ell},
\]
where $\sigma^2_{\ell}=\Var (f( X^{\ell,m^{\ell}}_{T})-f(
X^{\ell,m^{\ell-1}}_{T}) )$.
Assuming that the diffusion coefficients of $ X $ and the function $ f
$ are Lipschitz continuous, then it is easy to check, using properties
of the Euler scheme that
\[
\sigma^2\leq c_2\sum_{\ell=0}^LN^{-1}_{\ell}m^{-\ell}
\]
for some positive constant $c_2$ (see Proposition \ref{p1} for more
details). Giles \cite{Gil} uses this computation in order to find the
optimal choice of the multilevel Monte Carlo parameters. More precisely,
to obtain a desired root mean squared error (RMSE), say of order
$1/n^{\alpha}$,
for the multilevel estimator, Giles \cite{Gil} uses the above
computation on $\sigma^2$
to minimize the total complexity of the algorithm. It turns out that
the optimal choice is obtained
for (see Theorem~3.1 of \cite{Gil})
%e1 #&#
%
\begin{equation}
\label{GilesParam} N_{\ell}=2c_2n^{2\alpha} \biggl(
\frac{\log n}{\log m}+1 \biggr)\frac
{T}{m^{\ell}}\qquad\mbox{for } \ell\in\{0,\ldots,L\}
\mbox{ and }  L=\frac{\log n}{\log m}.
\end{equation}
Hence, for an error $\varepsilon_n=1/n^{\alpha}$, this optimal choice
leads to a complexity
for the multilevel Monte Carlo Euler method proportional to $n^{2\alpha
}(\log n)^2=\varepsilon_n^{-2}(\log\varepsilon_n)^2$.
Interesting numerical tests, comparing three methods (crude Monte
Carlo, statistical Romberg and the multilevel Monte Carlo), were
processed in Korn, Korn and Kroisandt \cite{Kor}.

In the present paper, we focus on central limit theorems for the
inferred error; a question which has not been addressed in previous
research.
To do so, we use techniques adapted to this setting, based on a central
limit theorem for triangular array
(see Theorem \ref{lindeberg}) together with Toeplitz lemma. It is worth
to note that our approach improves techniques developed by Kebaier
\cite{Keb}
in his study of the statistical Romberg method (see Remark \ref
{rem-keb} for more details).
Hence, our main result is a Lindeberg--Feller central limit theorem
for the multilevel Monte Carlo Euler algorithm (see Theorem \ref
{CLTSR}). Further, this allows us to prove a Berry--Esseen-type bound on
our central limit theorem.

In order to show this central limit theorem, we first prove a
stable law convergence theorem, for the Euler scheme error on two
consecutive levels $m^{\ell-1}$~and~$m^{\ell}$,
of the type obtained in Jacod and Protter \cite{Jacod-Protter}.
Indeed, we prove the following functional result
(see Theorem \ref{th-acc}):
\[
\sqrt{\frac{m^{\ell}}{(m-1)T}}\bigl(X^{\ell,m^{\ell}}-X^{\ell,m^{\ell-1}}\bigr)
\Rsta  U\qquad\mbox{as } \ell\to\infty,
\]
where $U$ is the same limit process given in Theorem 3.2 of Jacod and
Protter \cite{Jacod-Protter}.
Our result uses standard tools developed in their paper but it cannot
be deduced without a specific and laborious study.
Further, their result,
namely
\[
\sqrt{\frac{m^{\ell}}{T}}\bigl(X^{\ell,m^{\ell}}-X\bigr)\Rsta U\qquad\mbox{as } \ell\to\infty,
\]
is
neither sufficient nor appropriate to prove our Theorem \ref{CLTSR},
since the multilevel Monte Carlo Euler method involves the error process
$X^{\ell,m^{\ell}}-X^{\ell,m^{\ell-1}}$ rather than $X^{\ell,m^{\ell}}-X$.

Thanks to Theorem \ref{CLTSR}, we obtain a precise description for the
choice of the parameters to run the multilevel
Monte Carlo Euler method. Afterward, by a complexity analysis we obtain
the optimal choice for
the multilevel Monte Carlo Euler method. It turns out that for a total
error of order $\varepsilon_n=1/n^{\alpha}$ the optimal
parameters are given by
%e2 #&#
%
\begin{equation}
\label{BAKParam} N_{\ell}=\frac{(m -1)T}{m^\ell\log m}n^{2\alpha}\log n\qquad\mbox{for } \ell\in\{0,\ldots,L\}\mbox{ and } L=\frac{\log n}{\log m}.
\end{equation}
This leads us to a complexity proportional to $n^{2\alpha}(\log
n)^2=\varepsilon_n^{-2}(\log\varepsilon_n)^2$ which is the same order
obtained by Giles \cite{Gil}. By comparing relations
(\ref{GilesParam}) and (\ref{BAKParam}), we note that our optimal
sequence of sample sizes $(N_{\ell})_{0\leq\ell\leq L}$ does not
depend on
any given constant, since our approach is based on proving a central
limit theorem and not on obtaining an upper bound for the variance of
the algorithm. However, some numerical tests comparing the runtime with
respect to the root mean square error, show that we are in line with
the original work of Giles \cite{Gil}.
Nevertheless, the major advantage of our central limit theorem is that
it fills the gap in the literature for the multilevel Monte Carlo Euler
method and allows to construct a more accurate confidence interval
compared to the one obtained using Chebyshev's inequality.
All these results are stated and proved in Section~\ref{main}. The next
section is devoted to recall some useful stochastic limit theorems and
to introduce our notation.

%s2 #&#
\section{General framework}\label{sec2}
%s2.1 #&#
\subsection{Preliminaries}
Let $(X_n)$ be a sequence of random variables with values in a Polish
space $E$ defined on a probability
space $(\Omega,\mathcal F,\mathbb P)$. Let $(\widetilde\Omega,\widetilde{\mathcal F},\widetilde{\mathbb P})$ be an extension
of $(\Omega,\mathcal F,\mathbb P)$, and let $X$ be an $E$-valued random
variable on the extension.
We say that $(X_n)$ converges in law to $X$ stably and write
$X_n \Rightarrow^{\stably}X$, if
\[
\mathbb E\bigl(Uh(X_n)\bigr)\rightarrow\widetilde{\mathbb E}\bigl(Uh(X)
\bigr)
\]
for all $h\dvtx E\rightarrow\mathbb R$ bounded continuous and all bounded
random variable $U$ on $(\Omega,\mathcal F)$.
This convergence is obviously stronger than
convergence in law
that we will denote here by ``$\Rightarrow$.'' According to Section~2
of Jacod \cite{Jacod} and Lemma 2.1 of Jacod and Protter
\cite{Jacod-Protter}, we have the following result.

%le1 #&#
%
\begin{lemma}\label{lemma}
Let $V_n$ and $V$ be defined on $(\Omega,\mathcal F)$ with values in
another metric space $E'$.
\[
\mbox{If } V_n\stackrel{\mathbb P} {\rightarrow}V, X_n
\Rsta X\qquad\mbox{then } (V_n,X_n)
\Rsta (V,X).
\]
Conversely, if $(V,X_n)\Rightarrow(V,X)$ and $V$ generates the $\sigma
$-field ${\mathcal F}$, we can realize this limit
as $(V,X)$ with $X$ defined on an extension of $(\Omega,\mathcal
F,\mathbb P
)$ and $X_n \Rightarrow^{\stably} X$.
\end{lemma}

Now, we recall a result on the convergence of stochastic integrals
formulated from Theorem 2.3 in Jacod and Protter \cite{Jacod-Protter}.
This is a simplified version but it is sufficient for our study. Let
$X^n=(X^{n,i})_{1\leq i \leq d}$
be a sequence of $\mathbb R^d$-valued continuous semimartingales with
the decomposition
\[
X^{n,i}_t=X_0^{n,i}+A_t^{n,i}+M_t^{n,i},
\qquad0\leq t \leq T,
\]
where, for each $n\in\mathbb N$ and $1\leq i \leq d$, $A^{n,i}$ is a
predictable process with finite variation, null at $0$ and $M^{n,i}$ is
a martingale null at $0$.
%th1 #&#

\begin{theorem}
\label{th-conv-int}
Assume that the sequence $(X^n)$ is such that
\[
\label{eq-*} \bigl\langle M^{n,i}\bigr\rangle_T + \int
_{0}^{T}\bigl\llvert \,dA^{n,i}_s
\bigr\rrvert
\]
is tight. Let $H^n$ and $H$ be a sequence of adapted, right-continuous
and left-hand side limited processes all defined on the same filtered
probability space.
If $(H^n, X^n)\Rightarrow(H,X)$ then $X$ is a semimartingale with
respect to the filtration generated
by the limit process $(H,X)$, and we have
$
(H^n,X^n,\int H^nd X^n)\Rightarrow(H,X, \int H \,dX)$.
\end{theorem}

We recall also the following Lindeberg--Feller central limit theorem
that will be used in the sequel (see, e.g., Theorems 7.2~and~7.3 in~\cite{Bil}).

\begin{theorem}[(Central limit theorem for triangular array)]\label{lindeberg}
Let $(k_{n})_{n\in\mathbb N}$ be a sequence such that $k_{n}
\rightarrow
\infty$  as $n \rightarrow\infty$.
For each $n$,
let $X_{n,1}, \ldots, X_{n,k_n}$ be $k_n$ independent random variables
with finite variance such that $ \mathbb E(X_{n,k})=0$ for all $k\in\{
1,\ldots,k_n\}$. Suppose that the following conditions hold:
\begin{longlist}[(A3)]
\item[(A1)]
$\lim_{n\rightarrow\infty} \sum_{k=1}^{k_{n}} \mathbb{E} |X_{n,k}|^{2}
= \sigma^2, \sigma>0$.
\item[(A2)] Lindeberg's condition: for all $\varepsilon>0$,
$\lim_{n\rightarrow\infty}\sum_{k=1}^{k_{n}} \mathbb{E}  (|X_{n,k}
|^{2}\times\break 1_{\{|X_{n,k}|>\varepsilon\}}  )= 0$.
Then
\[
\sum_{k=1}^{k_{n}}X_{n,k}
\Rightarrow\mathcal{N}\bigl(0,\sigma^2\bigr)\qquad\mbox{as } n
\rightarrow\infty.
\]
\end{longlist}
Moreover, if the $X_{n,k}$ have moments of order $p>2$, then the
Lindeberg's condition can be obtained by the following one:
\begin{longlist}[(A3)]
\item[(A3)] Lyapunov's condition:
$ \lim_{n\rightarrow\infty} \sum_{k=1}^{k_{n}} \mathbb{E} |X_{n,k}|^{p}
= 0$.
\end{longlist}
\end{theorem}

%s2.2 #&#
\subsection{The Euler scheme}

Let $X:=  (X_t )_{\scriptscriptstyle{0\leq t\leq
T}}$ be the process with values in $\mathbb R^d$, solution to
%e3 #&#
%
\begin{equation}
\label{1} dX_{t} = b(X_{t})\,dt + {{\sigma(X_{t})}}\,dW_{t},\qquad
X_0=x\in\mathbb R^d,
\end{equation}
where $ W=(W^1,\dots,W^q)$ is a $q$-dimensional Brownian
motion on some given filtered probability space
$\mathcal B=(\Omega,\mathcal F,(\mathcal F_t)_{t\geq0},\mathbb P)$ with
$(\mathcal F_t)_{t\geq0}$ is the standard \mbox{Brownian} filtration, $b$ and
$\sigma$ are, respectively,
$\mathbb R^d$ and $\mathbb R^{d\times q}$ valued functions.
We consider the continuous Euler approximation $X^n$ with step $\delta
=T/n$ given by
\[
dX^n_t=b(X_{\eta_n(t)})\,dt+\sigma(X_{\eta_n(t)})\,dW_t,\qquad
\eta _n(t)=[t/\delta]\delta.
\]
It is well known that under the global Lipschitz condition
{\renewcommand{\theequation}{$\mathcal H_{b,\sigma}$}
\begin{eqnarray}\label{Hbsigma}
\qquad\exists C_T>0,\mbox{ such that, }
\bigl|b(x)-b(y)\bigr|+ \bigl|\sigma(x)-\sigma(y)\bigr|\leq C_T|y-x|,
\nonumber\\[-10pt]\\[-10pt]
\eqntext{x,y\in\mathbb R^d,}
\end{eqnarray}}%
the Euler scheme satisfies the following property (see, e.g., Bouleau
and L\'{e}pingle~\cite{BouLep}):
{\renewcommand{\theequation}{${\mathcal P}$}
\begin{eqnarray}\label{Pbold}
\forall p\geq1,\qquad &&\sup_{0 \leq t
\leq T}|X_t|,\qquad
\sup_{0 \leq t \leq T}\bigl|X^n_t\bigr|\in L^p
\quad\mbox{and}
\nonumber\\[-8pt]\\[-8pt]
&&\mathbb {E} \Bigl[{\sup_{0 \leq t \leq T}}\bigl|
X_t - {X}^n_t\bigr|^{p} \Bigr]
\leq\frac{K_p(T)}{n^{p/2}},\qquad K_p(T)>0.\nonumber
\end{eqnarray}}%
Note that according to Theorem 3.1 of Jacod and Protter \cite
{Jacod-Protter}, under the weaker condition
{\renewcommand{\theequation}{$\mathcal{\widetilde H}_{b,\sigma}$}
\begin{equation}\label{tildeHbsigma}
b\mbox{ and } \sigma \mbox{ are locally Lipschitz with linear growth},
\end{equation}}%
we have only the uniform convergence in probability, namely the property
{\renewcommand{\theequation}{$\widetilde{{\mathcal{P}}}$} %\bolds{\mathcal{\widetilde P
\begin{equation}\label{cPtilde}
\sup_{0 \leq t \leq T}\bigl|
X_t - {X}^n_t\bigr|\stackrel{\mathbb P} {\rightarrow}0.
\end{equation}}\setcounter{equation}{3}%
Following the notation of Jacod and Protter \cite{Jacod-Protter}, we
rewrite diffusion (\ref{1}) as follows:
\[
dX_t=\varphi(X_t)\,dY_t=\sum
_{j=0}^{q}\varphi_j(X_t)\,dY_t^j,
\]
where $\varphi_j$ is the $j$th column of the matrix $\sigma$, for
$1\leq j \leq q$, $\varphi_0=b$ and $Y_t:=(t,W^1_t,\ldots,W^q_t)'$.
Then the continuous Euler approximation $X^n$ with time step $\delta
=T/n$ becomes
%e4 #&#
%
\begin{equation}
\label{forme-euler} dX^n_t=\varphi\bigl(X^n_{\eta_n(t)}
\bigr)\,dY_t=\sum_{j=0}^{q}
\varphi _j\bigl(X^n_{\eta
_n(t)} \bigr)\,dY_t^j,\qquad
\eta_n(t)=[t/\delta]\delta.
\end{equation}
%
%s3 #&#
\section{The multilevel Monte Carlo Euler method}\label{main}
Let $(X^{m^{\ell}}_t)_{0\leq t\leq T}$ denotes the Euler scheme with
time step $m^{-\ell}T$ for $\ell\in\{0,\ldots,L\}$,
where $L=\log n/\break \log m$. Noting that
%e5 #&#
%
\begin{equation}
\label{eq-biais} \mathbb Ef\bigl( X^n_T\bigr)=\mathbb Ef
\bigl( X^1_T\bigr)+\sum_{\ell=1}^{L}
\mathbb E \bigl(f\bigl( X^{m^{\ell
}}_T\bigr)-f\bigl(
X^{m^{\ell-1}}_T\bigr) \bigr),
\end{equation}
the multilevel method is to estimate independently by the Monte Carlo
method each of the expectations on the right-hand side of the above relation.
Hence, we approximate $\mathbb Ef( X^n_T) $ by
%e6 #&#
%
\begin{equation}
\label{eq-estimateur} Q_n=\frac{1}{N_0}\sum
_{k=1}^{N_0}f\bigl( X^{1}_{T,k}
\bigr)+\sum_{\ell=1}^{L}\frac{1}{N_\ell}\sum
_{k=1}^{N_\ell} \bigl(f\bigl(
X^{\ell,m^{\ell
}}_{T,k}\bigr)-f\bigl( X^{\ell,m^{\ell-1}}_{T,k}
\bigr) \bigr).
\end{equation}
Here, it is important to point out that all these $L+1$ Monte
Carlo estimators have to be based on different, independent samples.
For each $\ell\in\{1,\ldots,L\}$ the samples $( X^{\ell,m^{\ell
}}_{T,k}, X^{\ell,m^{\ell-1}}_{T,k})_{1\leq k\leq N_{\ell}}$ are independent
copies of $(X^{\ell,m^{\ell}}_{T},X^{\ell,m^{\ell-1}}_{T})$ whose
components denote the Euler schemes with time steps $m^{-\ell}T$ and
$m^{-(\ell-1)}T$ and simulated with the same Brownian path. Concerning
the first empirical mean, the samples $(X^{1}_{T,k})_{1\leq k\leq N_{0}}$
are independent copies of $X^{1}_{T}$.
The following result gives us a first description of the asymptotic
behavior of the variance in the multilevel Monte Carlo Euler
method.

%pr1 #&#
%
\begin{proposition}\label{p1}
Assume that $b$ and $\sigma$ satisfy condition (\ref{Hbsigma}).
For a Lipschitz continuous function $ f\dvtx  \mathbb
R^d\longrightarrow\mathbb R$, we have
%e7 #&#
%
\begin{equation}
\label{ordre} \Var(Q_n) = O \Biggl(\sum_{\ell=0}^{L}N_{\ell}^{-1}m^{-\ell}
\Biggr).
\end{equation}
\end{proposition}
%
%pa3.subsection.subsubsection.1 #&#

\begin{pf}
We have
\begin{eqnarray*}
\Var(Q_n) &=& N_0^{-1} \Var \bigl(f\bigl(
X^1_T\bigr) \bigr) +\sum_{\ell=1}^{L}
N_{\ell
}^{-1} \Var \bigl(f\bigl( X^{\ell,m^{\ell}}_T
\bigr)-f\bigl( X^{\ell,m^{\ell
-1}}_T\bigr) \bigr)
\\
&\leq& N_0^{-1} \Var \bigl(f\bigl( X^1_T
\bigr) \bigr)
\\
&&{} +2\sum_{\ell=1}^{L}
N_{\ell
}^{-1} \bigl(\Var \bigl(f\bigl( X^{m^{\ell}}_T
\bigr)-f( X_T)\bigr)
+\Var \bigl(f\bigl( X^{m^{\ell
-1}}_T
\bigr)-f( X_T)\bigr) \bigr)
\\
&\leq& N_0^{-1} \Var \bigl(f\bigl( X^1_T
\bigr) \bigr)
\\
&&{} +2 [f]_{\mathrm{lip}} \sum_{\ell=1}^{L}
N_{\ell}^{-1}\mathbb E \Bigl[ \sup_{0 \leq t \leq T}
\bigl\llvert X^{m^{\ell
}}_t-X_t\bigr\rrvert ^2 + \sup_{0\leq t \leq T}\bigl\llvert X^{m^{\ell-1}}_t-X_t
\bigr\rrvert ^2 \Bigr],
\end{eqnarray*}
where\vspace*{2pt} $ [f]_{\mathrm{lip}}:=\sup_{u\neq v}\frac{|f(u)-f(v)|}{|u-v|}$. We complete
the proof by using property (\ref{Pbold}) on the strong
convergence of the Euler scheme.
\end{pf}

Inequality (\ref{ordre}) indicates the dependence of the variance of $Q_n$
on the choice of the parameters $N_{0},\dots,N_{L}$.
This variance can be smaller than the variance of $f( X^n_T)$, so that $Q_n$
appears as a good candidate for the variance reduction.

The main result of this section is
a Lindeberg--Feller central limit theorem (see Theorem \ref{CLTSR} below).
In order to prove this result, we need to prove first a new stable law
convergence theorem for the Euler scheme error adapted to the setting
of multilevel Monte Carlo algorithm. This is crucial and is the aim
of the following subsection.

%s3.1 #&#
\subsection{Stable convergence}\label{section2bis}
In what follows, we prove a
stable law convergence theorem, for the Euler scheme error on two
consecutive levels $m^{\ell-1}$ and $m^{\ell}$,
of the type obtained in Jacod and Protter \cite{Jacod-Protter}.
Our result in Theorem \ref{th-acc} below is an innovative contribution
on the Euler scheme error that is different and more tricky
than the original work by Jacod and Protter \cite{Jacod-Protter} since
it involves the error process
$X^{\ell,m^{\ell}}-X^{\ell,m^{\ell-1}}$ rather than $X^{\ell,m^{\ell}}-X$.
Note that the study of the error $X^{\ell,m^{\ell}}-X^{\ell,m^{\ell
-1}}$ as
$\ell\rightarrow\infty$ can be reduced to the study of the error $X^{m
n}-X^{n}$ as $n\rightarrow\infty$ where $X^{m n}$ and $X^{n}$ stand for
the Euler schemes with time steps $T/(mn)$ and $T/n$ constructed on
the same Brownian path.

%th2 #&#
%
\begin{theorem}\label{th-acc}
Assume that $b$ and $\sigma$ are $\mathcal C^1$ with linear growth
then the following result holds:
\[
\mbox{For all } m\in\mathbb N\setminus\{0,1\}\qquad\sqrt {
\frac{m
n}{(m-1)T}}\bigl(X^{m n}-X^{n}\bigr)
\Rsta  U\qquad\mbox{as }n\rightarrow\infty,
\]
with $(U_t)_{0\leq t\leq T}$ the $d$-dimensional process satisfying
%e8 #&#
%
\begin{equation}
\label{U1} U_t=\frac{1}{\sqrt{2}}\sum
^{q}_{i,j=1}Z_t\int_0^t
H^{i,j}_s \,dB^{ij}_s, \qquad t
\in[0,T],
\end{equation}
where
%e9 #&#
%
\begin{equation}
\label{eq-defH} \quad H^{i,j}_s=(Z_s)^{-1}
\dot\varphi_{s,j}\bar\varphi_{s,i}\qquad\mbox{with } \dot
\varphi_{s,j}:=\nabla\varphi_j (X_{s})\mbox{ and } \bar \varphi _{s,i}:=\varphi_i (X_{s}),
\end{equation}
and $(Z_t)_{0\leq t\leq T}$ is the $\mathbb R^{d\times d}$ valued process
solution of the linear equation
\[
Z_t=I_d+ \sum^{q}_{j=0}
\int_0^t \dot\varphi_{s,j}
\,dY^j_s Z_s,\qquad t\in[0,T].
\]
Here,\vspace*{1pt} $\nabla\varphi_{j}$ is a $d\times d$ matrix with $(\nabla
\varphi
_{j})_{ik}$ is the partial derivative of $\varphi_{ij}$ with respect to
the $k$th coordinate,
and $(B^{ij})_{1\leq i,j\leq q}$ is a standard $q^2$-dimensional
Brownian motion independent of $W$. This process is defined on an extension
$(\widetilde\Omega,\widetilde{\mathcal F},(\widetilde{\mathcal F}_t)_{t\geq
0},\widetilde{\mathbb P})$ of the space $(\Omega,\mathcal F,(\mathcal
F_t)_{t\geq
0},\mathbb P)$.
\end{theorem}

Note that by letting formally $m$ tend to infinity, we recover the
Jacod and Protter's result \cite{Jacod-Protter}.
%pa3.1.subsubsection.1 #&#

\begin{pf*}{Proof of Theorem \ref{th-acc}}
Consider the error process $U^{m n,n}=(U^{m n,n}_t)_{0\leq t\leq T}$,
defined by
\[
U^{m n,n}_t:=X^{m n}_t-X^n_t,
\qquad t\in[0,T].
\]
Combining relation (\ref{forme-euler}), for both processes $X^{m n}$
and $X^{n}$, together with a Taylor expansion
\[
dU^{m n,n}_t =\sum_{j=0}^{q}
\dot\varphi^{n}_{t,j}\bigl(X^{m n}_{\eta_{m
n}(t)}-X^n_{\eta
_n(t)}
\bigr) \,dY^j_t,
\]
where
$\dot\varphi^n_{t,j}$ is the $d\times d$ matrix whose $i$th row is the
gradient of the real-valued function $\varphi_{ij}$ at a point between
$X^n_{\eta_n(t)}$ and $X^{m n}_{\eta_{m n}(t)}$.
Therefore, the equation satisfied by $U^n$ can be written as
\[
\label{D1} U^{m n,n}_t=\int_0^t
\sum_{j=0}^{q}\dot\varphi^n_{s,j}
U^{m n,n}_s \,dY^j_s+G^{m n,n}_t,
\]
with
\[
\label{D2} G^{m n,n}_t=\int_0^t
\sum_{j=0}^{q}\dot\varphi^n_{s,j}
\bigl(X^n_s-X^n_{\eta_n(s)}\bigr)
\,dY^j_s -\int_0^t
\sum_{j=0}^{q}\dot\varphi^n_{s,j}
\bigl(X^{m n}_s-X^{m n}_{\eta_{m n}(s)}\bigr)
\,dY^j_s.
\]
In the following, let
$(Z^{m n,n}_t)_{0\leq t\leq T}$ be the $\mathbb{R}^{d\times d}$ valued
solution of
\[
Z^{m n,n}_t=I_d+\int_0^t
\Biggl(\sum_{j=0}^{q}\dot
\varphi^n_{s,j} \,dY^j_s \Biggr)
Z^{m n,n}_s.
\]
Theorem 48, page 326 in \cite{Protter}, ensures existence of the
process $ ((Z^{m n,n}_t)^{-1} )_{0\leq t\leq T}$ defined as
the solution of
\[
\bigl(Z^{m n,n}_t\bigr)^{-1}=I_d+\int
_0^t\bigl(Z^{m n,n}_s
\bigr)^{-1}\sum_{j=1}^{q}\bigl(
\dot \varphi^n_{s,j}\bigr)^2\,ds-\int
_0^t\bigl(Z^{m n,n}_s
\bigr)^{-1}\sum_{j=0}^{q}\dot
\varphi^n_{s,j}\,dY_s^j.
\]
Thanks to Theorem 56, page 333 in the same reference \cite{Protter},
we get
\begin{eqnarray*}
U^{m n,n}_t &=&Z^{m n,n}_t \Biggl\{ \int
_0^t \bigl(Z^{m n,n}_s
\bigr)^{-1}\,dG^{m n,n}_s
\\
&&\hspace*{32pt}{} -\int_0^t
\bigl(Z^{m n,n}_s\bigr)^{-1}\sum
_{j=1}^q\bigl(\dot\varphi^n_{s,j}
\bigr)^2 \bigl(X^n_s- X^n_{\eta
_n(s)}
\bigr) \,ds
\\
&&\hspace*{32pt}{}+\int_0^t \bigl(Z^{m n,n}_s
\bigr)^{-1}\sum_{j=1}^q\bigl(
\dot\varphi^n_{s,j}\bigr)^2
\bigl(X^{m n}_s-X^{m n}_{\eta_{m n}(s)}\bigr) \,ds
\Biggr\}.
\end{eqnarray*}
Since the increments of the Euler scheme satisfy
\[
X_s^n-X^n_{\eta_n(s)}=\sum
_{i=0}^{q}\bar\varphi ^n_{s,i}
\bigl(Y^i_s-Y^i_{\eta
_n(s)}\bigr)
\]
and
\[
X^{m n}_s-X^{m n}_{\eta_{m n}(s)} =
\sum_{i=0}^{q}\bar\varphi^{m
n}_{s,i}
\bigl(Y^i_s-Y^i_{\eta_{m n}(s)}\bigr),
\]
with $\bar\varphi^n_{s,i}=\varphi_i(X^n_{\eta_n(s)})$ and $\bar
\varphi^{m n}_{s,i}=\varphi_i( X^{m n}_{\eta_{m n}(s)} )$, it is easy
to check that
\begin{eqnarray}
\label{eq1erreur} \qquad U^{m n,n}_t &=& \sum
_{i,j=1}^{q}Z^{m n,n}_t \int
_0^t H^{i,j,m n,n}_s
\bigl(Y^i_s-Y^i_{\eta
_n(s)}\bigr)
\,dY^j_s+ R_{t,1}^{m n,n}+R_{t,2}^{m n,n}
\nonumber\\[-8pt]\\[-8pt]
&&{}-\sum_{i,j=1}^{q} Z^{m n,n}_t
\int_0^t \widetilde H^{i,j,m n,n}_s
\bigl(Y^i_s-Y^i_{\eta_{m n}(s)}\bigr)
\,dY^j_s -\widetilde R_{t,1}^{m n,n} -\widetilde
R_{t,2}^{m n,n}\nonumber
\end{eqnarray}
with
\begin{eqnarray*}
R_{t,1}^{m n,n}&=&\sum
_{i=0}^{q}Z^{m n,n}_t\int
_0^t K^{i,m
n,n}_s
\bigl(Y^i_s-Y^i_{\eta_n(s)}\bigr) \,ds,
\\
R_{t,2}^{m n,n}&=&\sum
_{j=1}^{q}Z^{m n,n}_t\int
_0^t H^{0,j,m n,n}_s \bigl(s-
\eta_n(s)\bigr) \,dY^j_s,
\end{eqnarray*}
and
\begin{eqnarray*}
\widetilde R_{t,1}^{m n,n} &=&\sum_{i=0}^{q}
Z^{m n,n}_t\int_0^t \widetilde
K^{i,m n,n}_s\bigl(Y^i_s-Y^i_{\eta_{m n}(s)}
\bigr)\,ds,
\\
\widetilde R_{t,2}^{m n,n}&=&\sum
_{j=1}^{q} Z^{m n,n}_t\int
_0^t \widetilde H^{0,j,m n,n}_s
\bigl(s-\eta_{m n}(s)\bigr)\,dY^j_s,
\end{eqnarray*}
where, for $(i,j)\in\{0,\ldots,q\}\times\{1,\ldots,q\}$,
\begin{eqnarray*}
K^{i,m n,n}_s&=&\bigl(Z^{m n,n}_s
\bigr)^{-1} \Biggl(\dot\varphi^n_{s,0}\bar \varphi
^n_{s,i} -\sum_{j=1}^q
\bigl(\dot\varphi^n_{s,j}\bigr)^2 \bar
\varphi^n_{s,i} \Biggr),
\\
H^{i,j,m n,n}_s&=&
\bigl(Z^{m n,n}_s\bigr)^{-1} \dot
\varphi^n_{s,j}\bar\varphi ^n_{s,i},
\end{eqnarray*}
and
\begin{eqnarray*}
\widetilde K^{i,m n,n}_s&=&\bigl(Z^{m n,n}_s
\bigr)^{-1} \Biggl(\dot\varphi^n_{s,0}\bar
\varphi^{m n}_{s,i}-\sum_{j=1}^q
\bigl(\dot\varphi^n_{s,j}\bigr)^2 \bar \varphi
^{m n}_{s,i} \Biggr),
\\
\widetilde H^{i,j,m n,n}_s&=&
\bigl(Z^{m n,n}_s\bigr)^{-1}\dot
\varphi^n_{s,j}\bar \varphi^{m n}_{s,i}.
\end{eqnarray*}
Now, let us introduce
\[
Z_t=I_d+\int_0^t\sum
_{j=0}^{q} \bigl(\dot\varphi_{s,j}
\,dY^j_s \bigr) Z_s\qquad\mbox{with }\dot
\varphi_{t,j}=\nabla\varphi_j(X_t).
\]
Moreover, $ ((Z_t)^{-1} )_{0\leq t\leq T}$ exists and
satisfies the following explicit linear stochastic differential equation:
\[
(Z_t)^{-1}=I_d+\int_0^t(Z_s)^{-1}
\sum_{j=1}^{q}(\dot\varphi
_{s,j})^2\,ds-\int_0^t(Z_s)^{-1}
\sum_{j=0}^{q}\dot\varphi_{s,j}\,dY_s^j.
\]
Thanks to the uniform convergence in probability of the Euler scheme
and according to Theorem 2.5 in Jacod and Protter \cite{Jacod-Protter},
we have
%e10 #&#
%
\begin{equation}
\label{eq-DoleansDade} \sup_{0\leq t\leq T}\bigl\llvert Z^{m n,n}_t-Z_t
\bigr\rrvert \stackrel{\mathbb P }\rightarrow0\quad\mbox{and}\quad \sup
_{0\leq t\leq T}\bigl\llvert {\bigl(Z^{m n,n}_t
\bigr)^{-1}}-{(Z_t)^{-1}}\bigr\rrvert \stackrel{
\mathbb P}\rightarrow0.
\end{equation}

Furthermore, in relation (\ref{eq1erreur}), one can replace,
respectively, $H^{i,j,m n,n}_s$ and $\widetilde H^{i,j,m n,n}_s$ by
their common limit $H^{i,j}_s$ given by relation (\ref{eq-defH}).
So that relation (\ref{eq1erreur}) becomes
%e11 #&#
%
\begin{equation}
\label{eq1biserreur} U^{m n,n}_t =\sum_{i,j=1}^{q}Z^{m n,n}_t
\int_0^t H^{i,j}_s
\bigl(Y^i_{\eta_{m
n}(s)}-Y^i_{\eta_n(s)}\bigr)
\,dY^j_s +R^{m n,n}_t,
\end{equation}
with
\[
R^{m n,n}_t=R_{t,1}^{m n,n}+R_{t,2}^{m n,n}+R_{t,3}^{m n,n}-
\widetilde R_{t,1}^{m n,n}- \widetilde R_{t,2}^{m n,n}-
\widetilde R_{t,3}^{m n,n},
\]
where $R_{t,i}^{m n,n}$ and $\widetilde R_{t,i}^{m n,n}$, $i\in\{1,2\}$,
are introduced by relation
(\ref{eq1erreur}) and
\begin{eqnarray*}
R_{t,3}^{m n,n}&=&\sum_{i,j=1}^{q}Z^{m n,n}_t
\int_0^t \bigl(H^{i,j,m
n,n}_s-H^{i,j}_s
\bigr) \bigl(Y^i_s-Y^i_{\eta_n(s)}
\bigr) \,dY^j_s,
\\
\widetilde R_{t,3}^{m n,n}&=&\sum_{i,j=1}^{q}Z^{m n,n}_t
\int_0^t \bigl(\widetilde H^{i,j,m n,n}_s-H^{i,j}_s
\bigr) \bigl(Y^i_s-Y^i_{\eta_{m n}(s)}
\bigr) \,dY^j_s.
\end{eqnarray*}
The remainder term process $R^{m n,n}$ vanishes with rate $\sqrt{n}$ in
probability. More precisely, we have the following convergence result.

%le2 #&#
%
\begin{lemma}\label{lem}
The rest term introduced in relation (\ref{eq1biserreur}) is such that
\[
\sup_{0\leq t \leq T} \bigl|\sqrt{n} R_{t}^{m n,n} \bigr|
\]
converges
to zero in probability as $n$ tends to infinity.
\end{lemma}

For the reader's convenience, the proof of this lemma is postponed to
the end of the current subsection.

The task is now to study the asymptotic behavior of the process given
by relation~(\ref{eq1biserreur})
\[
\sum_{i,j=1}^{q}\sqrt{n}Z^{m n,n}_t
\int_0^t H^{i,j}_s
\bigl(Y^i_{\eta_{m
n}(s)}-Y^i_{\eta_n(s)}\bigr)
\,dY^j_s.
\]
In order to study this process, we introduce the martingale process,
\[
\label{eq-martingale} M_t^{n,i,j}=\int_0^t
\bigl(Y^i_{\eta_{m n}(s)}-Y^i_{\eta_n(s)}\bigr)
\,dY^j_s,\qquad(i,j)\in\{1,\ldots,q\}^2,
\]
and we proceed to a preliminary calculus of the expectation of its
bracket.

Let $(i,j)$ and $(i',j') \in\{1,\ldots,q\}^2$, we have:
\begin{itemize}
\item for $j\neq j'$, the bracket $\langle
M^{n,i,j},M^{n,i',j'}\rangle=0$,
\item for $j=j'$ and $i\neq i'$, $\mathbb E\langle
M^{n,i,j},M^{n,i',j}\rangle=0$,
\item for $j=j'$ and $i=i'$, $\mathbb E\langle M^{n,i,j}\rangle
_t=\int_0^t(\eta_{m n}(s)-\eta_n(s)) \,ds$, $t\in[0,T]$ and we have
%e12 #&#
%
\begin{eqnarray}
\label{eqn-expectation} \mathbb E\bigl(\bigl\langle M^{n,i,j} \bigr
\rangle_t\bigr)&=&\int_{0}^{\eta_n(t)}\bigl(
\eta_{m
n}(s)-\eta_n(s)\bigr) \,ds+O\biggl(\frac{1}{n^2}
\biggr)
\nonumber
\\
&=&\sum_{\ell=0}^{m-1}\sum
_{k=0}^{[t/\delta]-1}\int_{(mk+\ell
)\delta
/m}^{(mk+\ell+1)\delta/m}
\bigl(\eta_{m n}(s)-\eta_n(s)\bigr) \,ds+O\biggl(
\frac
{1}{n^2}\biggr)
\nonumber
\\
&=&\sum_{\ell=0}^{m-1}\sum
_{k=0}^{[t/\delta]-1}\frac{\delta
^2}{m} \biggl(
\frac{mk+\ell}{m}-k \biggr)+O\biggl(\frac{1}{n^2}\biggr)
\\
&=& \frac{(m-1)\delta
^2}{2m}[t/\delta] +O\biggl(\frac{1}{n^2}\biggr)
\nonumber
\\
&=&\frac{(m -1)T}{2m n} t +O\biggl(\frac{1}{n^2}\biggr).\nonumber
\end{eqnarray}
\end{itemize}

Having disposed of this preliminary evaluations, we can now study the
stable convergence of $ (\sqrt{\frac{2m n}{(m -1)T}
}M^{n,i,j}
)_{1\leq i,j\leq q}$.
By virtue of Theorem 2.1 in \cite{Jacod}, we need to study the
asymptotic behavior of both brackets
$n\langle M^{n,i,j},M^{n,i',j'}\rangle_t$ and $\sqrt{n}\langle
M^{n,i,j},Y^{j'} \rangle_t$, for all
$t\in[0,T]$ and all $(i,j,i',j') \in\{1,\ldots,q\}^4$.
The case \mbox{$j\neq j'$} is obvious and we only proceed to prove that:
\begin{itemize}
\item for $j=j'$, $\sqrt{n}\langle M^{n,i,j},Y^{j} \rangle_t \mathop{\longrightarrow}\limits _{n \to\infty}^{{\mathbb P}}0$, for all
$t\in[0,T]$,\vspace*{1pt}
\item for $j=j'$ and $i\neq i'$,
$n\langle M^{n,i,j},M^{n,i',j}\rangle_t  {\mathop
{\longrightarrow
}\limits _{n \to\infty}^{{\mathbb P}}}0$, for all $t\in[0,T]$,\vspace*{1pt}
\item for $j=j'$ and $i =i'$, $n\langle M^{n,i,j}\rangle_t
{\mathop{\longrightarrow}\limits _{n \to\infty}^{{\mathbb P}}} \frac{(m
-1)T}{2m } t$, for all $t\in[0,T]$.
\end{itemize}
For the first point, we consider the $L^2$ convergence
\begin{eqnarray*}
\mathbb E\bigl\langle M^{n,i,j},Y^{j} \bigr
\rangle^2_t &=&\mathbb E \biggl(\int_{0}^{t}
\bigl(Y^i_{\eta_{m n}(s)}-Y^i_{\eta_n(s)}\bigr)\,ds
\biggr)^2
\\
&=&\int_{0}^{t}\!\!\int_{0}^{t}
\mathbb E \bigl(\bigl(Y^i_{\eta_{m
n}(s)}-Y^i_{\eta
_n(s)}
\bigr) \bigl(Y^i_{\eta_{m n}(u)}-Y^i_{\eta_n(u)}
\bigr) \bigr) \,ds\,du
\\
&=&2 \int_{0<s<u<t}g(s,u)\,ds \,du
\end{eqnarray*}
with
%e13 #&#
%
\begin{eqnarray}\label{eq-defg}
g(s,u)&=&\eta_{m n}(s)\wedge\eta_{m n}(u)-
\eta_{m n}(s)\wedge\eta_n(u)
\nonumber\\[-8pt]\\[-8pt]
&&{}  -\eta_n(s)\wedge
\eta_{m n}(u)+ \eta_n(s)\wedge\eta_n(u).\nonumber
\end{eqnarray}
It is worthy to note that
%e14 #&#
%
\begin{equation}
\label{eq-prop} \eta_n(s)\leq\eta_{m n}(s) \leq s \leq
\eta_n(u) \leq\eta_{m n}(u) \leq u\qquad\forall s\leq
\eta_n(u).
\end{equation}
Hence, $g(s,u)=0$, for $s\leq\eta_n(u)$, $g(s,u)=\eta_{m n}(s)-\eta
_{n}(s)$, for $\eta_n(u)< s <u$, and
\begin{eqnarray*}
\mathbb E \bigl\langle M^{n,i,j},Y^{j} \bigr
\rangle^2_t &=& 2\int_{0<\eta_n(u)<s<u<t} \bigl(
\eta_{m n}(s)-\eta_{n}(s) \bigr)\,ds \,du
\\
&\leq& 2\frac{T}{n}
\int_{0}^{t}\bigl(u-\eta_n(u)\bigr)\,du
\\
&\leq& 2\frac{T^2}{n^2}t.
\end{eqnarray*}
This yields the desired result. Concerning the second point, the $L^2$
norm is given by
\begin{eqnarray*}
\mathbb E\bigl\langle M^{n,i,j},M^{n,i',j} \bigr
\rangle^2_t&=&\mathbb E \biggl(\int_{0}^{t}
\bigl(Y^i_{\eta_{m n}(s)}-Y^i_{\eta_n(s)}\bigr)
\bigl(Y^{i'}_{\eta_{m
n}(s)}-Y^{i'}_{\eta_n(s)}\bigr)\,ds
\biggr)^2
\\
&=&\int_{0}^{t}\!\!\int_{0}^{t}
\bigl(\mathbb E \bigl(\bigl(Y^i_{\eta_{m
n}(s)}-Y^i_{\eta
_n(s)}
\bigr) \bigl(Y^i_{\eta_{m n}(u)}-Y^i_{\eta_n(u)}
\bigr) \bigr) \bigr)^2 \,ds\,du
\\
&=&2 \int_{0<s<u<t}g(s,u)^2\,ds \,du,
\end{eqnarray*}
with the same function $g$ given in relation (\ref{eq-defg}). Using
the properties of function $g$ developed above, we have in the same manner
\[
\mathbb E\bigl\langle M^{n,i,j},M^{n,i',j} \bigr
\rangle^2_t =2\int_{0<\eta
_n(u)<s<u<t}\bigl(
\eta_{m n}(s)-\eta_{n}(s)\bigr)^2\,ds \,du\leq2
\frac{T^3}{n^3}t,
\]
which proves our claim. For the last point, that is the essential one,
we use the development of $\mathbb E\langle M^{n,i,j}\rangle_t$ given by
relation (\ref{eqn-expectation})
to get
%e15 #&#
%
\begin{eqnarray}\label{eq-dev1} 
&& \mathbb E \biggl(n\bigl\langle M^{n,i,j} \bigr
\rangle_t-\frac{(m -1)T}{2m} t \biggr)^2
\nonumber\\[-8pt]\\[-8pt]\nonumber
&&\qquad =n^2
\mathbb E\bigl\langle M^{n,i,j} \bigr\rangle_t^2-
\frac{(m -1)^2T^2}{4m^2} t^2+O\biggl(\frac{1}{n}\biggr).
\end{eqnarray}
Otherwise, we have
%e16 #&#
%
\begin{eqnarray}
\label{eqn2-expectation} \mathbb E\bigl\langle M^{n,i,j} \bigr\rangle^2_t
&=&\mathbb E \biggl(\int_{0}^{t}
\bigl(Y^i_{\eta_{m
n}(s)}-Y^i_{\eta_n(s)}
\bigr)^2\,ds \biggr)^2
\nonumber
\\
&=&\int_{0}^{t}\!\!\int_{0}^{t}
\mathbb E \bigl(\bigl(Y^i_{\eta_{m
n}(s)}-Y^i_{\eta
_n(s)}
\bigr)^2\bigl(Y^i_{\eta_{m n}(u)}-Y^i_{\eta_n(u)}
\bigr)^2 \bigr) \,ds\,du
\\
&=&2 \int_{0<s<u<t}h(s,u)\,ds \,du\nonumber
\end{eqnarray}
with
%e17 #&#
%
\begin{equation}
\label{eq-defh} h(s,u)=\mathbb E \bigl(\bigl(Y^i_{\eta_{m n}(s)}-Y^i_{\eta
_n(s)}
\bigr)^2\bigl(Y^i_{\eta_{m
n}(u)}-Y^i_{\eta_n(u)}
\bigr)^2 \bigr).
\end{equation}
On one hand, for $s\leq\eta_n(u)$, using property (\ref{eq-prop})
together with the independence of the increments $Y^i_{\eta_{m
n}(s)}-Y^i_{\eta_n(s)}$ and $Y^i_{\eta_{m n}(u)}-Y^i_{\eta_n(u)}$, yields
\[
h(s,u)=\bigl(\eta_{m n}(s)-\eta_n(s)\bigr) \bigl(
\eta_{m n}(u)-\eta_n(u)\bigr).
\]
On the other hand, in relation
(\ref{eq-defh}) we use the Cauchy--Schwarz inequality to get
$h(s,u)=O(\frac{1}{n^2})$ and this yields
\[
\int_{0<\eta_n(u)<s<u<t}h(s,u)\,ds \,du=O\biggl(\frac{1}{n^3}\biggr).
\]
Now, noting that $(\eta_{m n}(s)-\eta_n(s))(\eta_{m n}(u)-\eta
_n(u))=O(\frac{1}{n^2})$, relation (\ref{eqn2-expectation}) becomes
\begin{eqnarray*}
\mathbb E \bigl(\bigl\langle M^{n,i,j} \bigr\rangle^2_t
\bigr)&=&2 \int_{0<s<u<t}\bigl(\eta _{m n}(s)-
\eta_n(s)\bigr) \bigl(\eta_{m n}(u)-\eta_n(u)
\bigr) \,ds \,du + O\biggl(\frac
{1}{n^3}\biggr)
\\
&=& \biggl(\int_{0}^{t}\bigl(
\eta_{m n}(s)-\eta_n(s)\bigr) \,ds \biggr)^2+O
\biggl(\frac
{1}{n^3}\biggr).
\end{eqnarray*}
Once again thanks to the development of $\mathbb E(\langle
M^{n,i,j}\rangle_t)$
given by relation (\ref{eqn-expectation}), we deduce that
%e18 #&#
%
\begin{equation}
\label{eq-dev2} \mathbb E\bigl\langle M^{n,i,j} \bigr\rangle^2_t
=\frac{(m -1)^2 T^2}{4m^2 n^2} t^2+O\biggl(\frac{1}{n^3}\biggr).
\end{equation}
By (\ref{eq-dev1}) and (\ref{eq-dev2}), we deduce the
convergence in $L^2$ of $n\langle M^{n,i,j} \rangle_t$ toward $\frac{(m
-1)T}{2m} t$.
By Theorem 2.1 in Jacod \cite{Jacod}, $ (\sqrt{\frac
{2 m
n}{(m -1)T} }M^{n,i,j} )_{1\leq i,j\leq q}$ converges in law
stably to a standard $q^2$-dimensional Brownian motion $(B^{ij})_{1\leq
i,j\leq q}$ independent of $W$. Consequently, by Lemma \ref{lemma} and
Theorem \ref{th-conv-int}, we obtain
\begin{eqnarray*}
&& \biggl(\sqrt{\frac{m n}{(m -1)T} }\int_0^t
H^{i,j}_s \bigl(Y^i_{\eta_{m
n}(s)}-Y^i_{\eta_n(s)}
\bigr)\,dY^j_s, t\geq0 \biggr)_{1\leq i,j\leq q}
\\
&&\qquad \Rsta  \biggl( \int_0^t
H^{i,j}_s \frac{dB^{ij}_s}{\sqrt{2}}, t\geq0 \biggr)_{1\leq i,j\leq q}.
\end{eqnarray*}
Finally, we complete the proof using relations (\ref{eq-DoleansDade}),
(\ref{eq1biserreur}), Lemma \ref{lem} and once again
Lemma \ref{lemma} to obtain
\[
\sqrt{\frac{m n}{(m -1)T} } U^{m n,n} \Rsta  U\qquad\mbox{where } U_t=\frac{1}{\sqrt {2}}\sum_{i,j=1}^{q}
Z_t \int_0^t H^{i,j}_s
\,dB^{ij}_s.
\]\upqed
\end{pf*}

%pa3.1.subsubsection.2 #&#

\begin{pf*}{Proof of Lemma \protect\ref{lem}}
At first, we prove the uniform convergence in probability toward zero
of the normalized rest terms
$\sqrt{n}R_{t,i}^{m n,n}$ for $i\in\{1,2\}$. The convergence of
$\sqrt {n}\widetilde R_{t,i}^{m n,n}$ $i\in\{1,2\}$
is a straightforward consequence of the previous one. The main part of
these rest terms
can be represented as integrals with respect to three types of
supermartingales that can be classified through the following three cases:
\begin{eqnarray*}
D_t^{n,0,0} &=& \sqrt{n}\int_0^t
\bigl(s-\eta_n(s)\bigr) \,ds,
\\
D_t^{n,i,0}&=&\sqrt {n}
\int_0^t \bigl(Y^i_s-Y^i_{\eta_n(s)}
\bigr) \,ds,
\\
M_t^{n,0,j}&=&\sqrt{n}\int_0^t
\bigl(s-\eta_n(s)\bigr) \,dY^j_s,
\end{eqnarray*}
where $(i,j) \in\{1,\ldots,q\}^2$ and $t\in[0,T]$. In the first case,
the supermartingale is deterministic of finite variation and its total
variation on the interval $[0,T]$ has the following expression:
\[
\int_{0}^{T}\bigl\llvert \,dD_t^{n,0,0}
\bigr\rrvert =\sqrt{n} \int_{0}^{T}\bigl(s-\eta
_n(s)\bigr) \,ds\leq\frac{T^2}{\sqrt{n}}.
\]
So, the process $D^{n,0,0}$ converges to $0$ and is tight.
In the second case, for $i\in\{1,\ldots,q\}$, the supermartingale is
also of finite variation and its total variation on the interval
$[0,T]$ has the following expression:
\[
\int_{0}^{T}\bigl\llvert \,dD_t^{n,i,0}
\bigr\rrvert =\sqrt{n} \int_{0}^{T}\bigl|Y^i_s-Y^i_{\eta_n(s)}\bigr|
\,ds.
\]
It is clear that $\sup_{n} \mathbb E (\int_{0}^{T}\llvert \,dD_s^{n,i,0}\rrvert  )<\infty$,
which ensures the tightness of the process $D^{n,i,0}$. Therefore, we
only need to establish the convergence of $D_t^{n,i,0}$ toward~$0$ in
$L^2(\Omega)$, for $t\in[0,T]$. In fact, we have
\[
\mathbb E \bigl(\bigl(D_t^{n,i,0}\bigr)^2
\bigr)=2n\int_{0<s<u<t} \mathbb E \bigl(\bigl(Y^{i}_s-Y^{i}_{\eta_n(s)}
\bigr) \bigl(Y^{i}_u-Y^{i}_{\eta_n(u)}
\bigr) \bigr) \,ds \,du.
\]
When $s\leq\eta_n(u)$, we have $\eta_n(s)\leq s \leq\eta_n(u) \leq
u$ and
by independence of the Brownian motion increments, we deduce that the
integrand term is equal to $0$.
Otherwise, when $s\geq\eta_n(u)$, we apply the Cauchy--Schwarz
inequality to get
\[
\mathbb E \bigl(\bigl(D_t^{n,i,0}\bigr)^2 \bigr)
\leq2T\int_{0}^{t}\bigl(u-\eta _n(u)
\bigr)\,du\leq 2\frac{T^2}{n}t.
\]
It follows from all these that $D^{n,i,0}\Rightarrow0$.
In the last case, for $j\in\{1,\ldots,q\}$, the process $M_t^{n,0,j}$
is a square integrable martingale and its bracket has the following expression:
\[
\bigl\langle M^{n,0,j} \bigr\rangle_T=n \int
_0^T \bigl(s-\eta_n(s)
\bigr)^2 \,ds\leq\frac
{T^3}{n}.
\]
It is clear that $\sup_{n} \mathbb E\langle M^{n,0,j} \rangle_T
<\infty$, so
we deduce the tightness of the process $\langle M^{n,0,j} \rangle$ and
the convergence $M^{n,0,j}\Rightarrow0$.

Now thanks to property (\hyperref[cPtilde]{$\widetilde{\mathcal{P}}$}) and relation (\ref
{eq-DoleansDade}), it is easy to check
that the integrand processes $K^{i,m n,n}_s$ and $H^{0,j,m n,n}_s$,
introduced in relation (\ref{eq1erreur}),
converge uniformly in probability to their respective limits
$
K^{i}_s=(Z_s)^{-1} (\dot\varphi_{s,0}\bar\varphi_{s,i} -\sum_{j=1}^q(\dot\varphi_{s,j})^2 \bar\varphi_{s,i}  )
$
and\vspace*{2pt}
$
H^{0,j}_s=(Z_s)^{-1}\dot\varphi_{s,j}\bar\varphi_{s,i}$,
where $\dot\varphi_{s,j}=\nabla\varphi_j (X_{s})$ and $\bar\varphi
_{s,i}=\varphi_i (X_{s})$.
Therefore, by Theorem \ref{th-conv-int} we deduce that the integral
processes given by
\[
\sqrt{n}\int_0^t K^{i,m n,n}_s
\bigl(Y^i_s-Y^i_{\eta_n(s)}\bigr) \,ds
\quad\mbox{and}\quad \sqrt{n}\int_0^t
H^{0,j,m n,n}_s \bigl(s-\eta_n(s)\bigr)
\,dY^j_s
\]
vanish.
Consequently, we conclude using relation (\ref{eq-DoleansDade}) that
$\sqrt{n}R_{i}^{m n,n}\Rightarrow0$ for $i\in\{1,2\}$.

We now proceed to prove that $R_{3}^{m n,n}\Rightarrow0$. The
convergence of
the process $\widetilde R_{3}^{m n,n}$ toward $0$ is obviously obtained
from the previous one.
The main part of this rest term can be represented as a stochastic
integral with respect to the martingale process given by
\[
N_t^{n,i,j}=\sqrt{n}\int_0^t
\bigl(Y^i_{s}-Y^i_{\eta_n(s)}\bigr)
\,dY^j_s,
\]
with $(i,j)\in\{1,\ldots,q\}\times\{1,\ldots,q\}$.
It was proven in Jacod and Protter \cite{Jacod-Protter} that
\[
\sqrt{ \frac{n}{T}} N^{n,i,j}\Rsta
\frac
{B^{ij}}{\sqrt{2}},
\]
where\vspace*{2pt} $(B^{ij})_{1\leq i,j\leq q}$ is a standard $q^2$-dimensional
Brownian motion defined on an extension probability space
$(\widetilde\Omega,\widetilde{\mathcal F},(\widetilde{\mathcal F}_t)_{t\geq
0},\widetilde{\mathbb P})$ of $(\Omega,\mathcal F,(\mathcal F_t)_{t\geq
0},\mathbb P)$,
which\vspace*{2pt} is independent of $W$.
Thanks\vspace*{2pt} to property (\hyperref[cPtilde]{$\widetilde{\mathcal{P}}$}) and relation (\ref
{eq-DoleansDade}), the integrand process $H^{i,j,m
n,n}-H^{i,j}\Rightarrow0$ and once again
by Theorem \ref{th-conv-int} we deduce that the integral processes
given by
\[
\sqrt{n}\int_0^t \bigl(H^{i,j,m n,n}_s-H^{i,j}_s
\bigr) \bigl(Y^i_s-Y^i_{\eta
_n(s)}
\bigr) \,dY^j_s
\]
vanish.
All this allows us to conclude using relation (\ref{eq-DoleansDade}).
\end{pf*}

%s3.2 #&#
\subsection{Central limit theorem}\label{CLT-Euler}
Let us recall that the multilevel Monte Carlo method uses information
from a sequence of computations with decreasing step sizes and
approximates the quantity $\mathbb Ef( X_T)$ by
\begin{eqnarray}
Q_n=\frac{1}{N_0}\sum_{k=1}^{N_0}f
\bigl( X^{1}_{T,k}\bigr)+\sum_{\ell=1}^{L}
\frac{1}{N_\ell} \sum_{k=1}^{N_\ell} \bigl(f
\bigl( X^{\ell,m^{\ell}}_{T,k}\bigr)-f\bigl( X^{\ell,m^{\ell
-1}}_{T,k}
\bigr) \bigr),\nonumber
\\
\eqntext{\displaystyle m\in\mathbb N\setminus\{0,1\}\mbox{ and }L=\frac{\log n}{\log m}.}
\end{eqnarray}
In the same way as in the case of a crude Monte Carlo estimation, let
us assume that the discretization error
\[
\varepsilon_n=\mathbb Ef\bigl(X^n_T\bigr)-
\mathbb Ef(X_T)
\]
is of order $1/n^{\alpha}$ for any $\alpha\in[1/2,1]$.
Taking advantage from the limit theorem proven in the above section, we
are now able to establish a central limit theorem of Lindeberg--Feller type
on the multilevel Monte Carlo Euler method. To do so, we introduce a
real sequence $(a_\ell)_{\ell\in\mathbb N}$ of positive terms such that
{\renewcommand{\theequation}{$\mathcal W$}
\begin{eqnarray}\label{con-weights-e}
\lim_{L\rightarrow
\infty}
\sum_{\ell=1}^{L}a_\ell=\infty
\quad\mbox{and}\quad \lim_{L\rightarrow\infty}\frac{1}{ (\sum_{\ell
=1}^{L}a_\ell )^{p/2}}\sum
_{\ell=1}^{L}a_\ell^{p/2} =0
\nonumber\\[-10pt]\\[-10pt]
\eqntext{\mbox{for } p>2}
\end{eqnarray}}\setcounter{equation}{19}%
and we assume that the sample size $N_{\ell}$ depends on the rest of
parameters by the relation
%e19 #&#
%
\begin{equation}
\label{eq-size} N_{\ell}=\frac{n^{2\alpha}(m -1)T}{m^\ell a_\ell}\sum
_{\ell
=1}^{L}a_\ell, \qquad\ell\in\{0,\ldots,L\}\mbox{ and } L=\frac{\log
n}{\log m}.
\end{equation}
We choose this form for $N_{\ell}$ because it is a generic form
allowing us a straightforward use of Toeplitz lemma that is a crucial
tool used in the proof of our central limit theorem. Indeed, property
(\ref{con-weights-e}) implies that if $(x_{\ell})_{\ell\geq1}$ is a
sequence converging
to $x\in\mathbb R$ as $\ell$ tends to infinity then
\[
\lim_{L\rightarrow+\infty} \frac{\sum_{\ell=1}^{L} a_{\ell}
x_{\ell
}}{\sum_{\ell=1}^{L}a_\ell} =x.
\]
In the sequel, we will denote by $\widetilde{\mathbb E}$, respectively,
$\widetilde \Var$
the expectation,\vspace*{1pt} respectively, the variance defined on the probability space
$(\widetilde\Omega,\widetilde{\mathcal F},\widetilde{\mathbb P})$ introduced in Theorem
\ref{th-acc}.
We can now state the central limit theorem under strengthened
conditions on the diffusion coefficients.

%th3 #&#
%
\begin{theorem}\label{CLTSR}Assume that $b$ and $\sigma$ are
$\mathcal
C^1$ functions satisfying the global Lipschitz condition (\ref{Hbsigma}).
Let f be a real-valued function satisfying
{\renewcommand{\theequation}{${\mathcal H}_f$}% ({\bf\mathcal H_f})
\begin{equation}\label{Hf}
\qquad \bigl|f(x)-f(y)\bigr|\leq C \bigl(1+|x|^p+|y|^p\bigr)|x-y|\qquad\mbox{for some $C,p>0$}.
\end{equation}}%
Assume $\mathbb P(X_T\notin\mathcal{D}_f)=0$, where $\mathcal
{D}_f:=\{
x\in\mathbb R^d; f $ is differentiable at $x\}$,
and that for some $\alpha\in[1/2,1]$ we have
{\renewcommand{\theequation}{${\mathcal H}_{\varepsilon_n}$}  %({\bf\mathcal H_{\varepsilon_n}})
\begin{equation}\label{Hen}
\lim_{n\rightarrow\infty}
n^{\alpha}\varepsilon_n=C_f(T,\alpha).
\end{equation}}\setcounter{equation}{20}%
Then, for the choice of $N_{\ell}$, $\ell\in\{0,1,\ldots,L\}$ given by
equation (\ref{eq-size}), we have
\[
n^{\alpha} \bigl(Q_{n} - \mathbb E \bigl( f(X_T)
\bigr) \bigr)\Rightarrow \mathcal N \bigl(C_f(T,\alpha),
\sigma^2 \bigr)
\]
with $\sigma^2=\widetilde \Var (\nabla f(X_T).U_T )$ and
$\mathcal
N (C_f(T,\alpha), \sigma^2  )$ denotes a normal distribution.
\end{theorem}

The global Lipschitz condition (\ref{Hbsigma}) seems to be
essential to establish our result, since it ensures property~(\ref{Pbold}). Otherwise, Hutzenthaler, Jentzen and Kloeden~\cite{HutJenKlo} prove that under weaker conditions on $b$ and $\sigma$
the multilevel Monte Carlo Euler method may diverges whereas the crude
Monte Carlo method converges.
%pa3.2.subsubsection.1 #&#

\begin{pf*}{Proof of Theorem \ref{CLTSR}}
To simplify our notation, we give the proof for \mbox{$\alpha=1$}, the case
$\alpha\in[1/2,1)$ is a straightforward deduction.
Combining relations (\ref{eq-biais}) and (\ref{eq-estimateur})
together, we get
\[
Q_{n} - \mathbb E \bigl(f(X_T) \bigr) =\widehat Q^1_n+ \widehat Q^2_n+
\varepsilon_n,
\]
where
\begin{eqnarray*}
\widehat Q^1_n &= &\frac{1}{N_0}\sum
_{k=1}^{N_0} \bigl(f\bigl( X^1_{T,k}
\bigr)-\mathbb E \bigl(f\bigl( X^1_T\bigr) \bigr) \bigr),
\\
\widehat Q^2_n &=&\sum_{\ell=1}^{L}
\frac{1}{N_\ell}\sum_{k=1}^{N_\ell
} \bigl(f
\bigl( X^{\ell,m^{\ell}}_{T,k}\bigr)-f\bigl( X^{\ell,m^{\ell-1}}_{T,k}
\bigr)-\mathbb E \bigl(f\bigl( X^{\ell,m^{\ell}}_T\bigr)-f\bigl(
X^{\ell,m^{\ell-1}}_T\bigr) \bigr) \bigr).
\end{eqnarray*}
Using assumption (\ref{Hen}), we obviously
obtain the term $C_f(T,\alpha)$ in the limit. Taking $N_0=\frac
{n^2(m-1)T}{a_0}\sum_{\ell=1}^{L}a_\ell$, we\vspace*{-1pt} can apply the classical
central limit theorem to~$\widehat Q^1_n$. Then we have
$
n \widehat Q^1_n \stackrel{\mathbb P}\rightarrow0$.
Finally,\vspace*{3pt} we have only to study the convergence of $n \widehat Q^2_n$ and we
will conclude by establishing
\[
n \widehat Q^2_n \Rightarrow\mathcal N \bigl(0,\widetilde \Var
\bigl(\nabla f(X_T).U_T \bigr) \bigr).
\]
To do so, we plan to use Theorem \ref{lindeberg} with the Lyapunov
condition and
we set
%e20 #&#
%
\begin{eqnarray}\label{Xnl-BerryESSen}
X_{n,\ell} &:=& \frac{n}{N_\ell}\sum
_{k=1}^{N_\ell
}Z^{m^{\ell
},m^{\ell-1}}_{T,k}\quad\mbox{and}
\nonumber\\[-8pt]\\[-8pt]
Z^{m^{\ell},m^{\ell
-1}}_{T,k}&:=& f\bigl( X^{\ell,m^{\ell}}_{T,k}
\bigr)-f\bigl( X^{\ell,m^{\ell-1}}_{T,k}\bigr)-\mathbb E \bigl(f\bigl(
X^{\ell,m^{\ell}}_{T,K}\bigr)-f\bigl( X^{\ell,m^{\ell-1}}_{T,k}
\bigr) \bigr).\nonumber
\end{eqnarray}
In other words, we will check the following conditions:
\begin{itemize}
\item $\lim_{n\rightarrow\infty} \sum_{\ell=1}^{L}\mathbb E(X_{n,\ell
})^2=\widetilde
\Var (\nabla f(X_T).U_T )$.

\item (Lyapunov condition) there exists $p>2$ such that
$
\lim_{n\rightarrow\infty} \sum_{\ell=1}^{L}\mathbb E\llvert X_{n,\ell}\rrvert ^{p}=0$.
\end{itemize}

For the first one, we have
\begin{eqnarray}
\label{TCL-LindebergFeller}
\sum_{\ell=1}^{L}\mathbb
E(X_{n,\ell})^2&=&\sum_{\ell
=1}^{L}\Var(X_{n,\ell})\nonumber
\\
&=& \sum_{\ell=1}^{L}\frac{n^2}{N_\ell}\Var
\bigl(Z^{m^{\ell},m^{\ell
-1}}_{T,1} \bigr)
\\
&=&\frac{1}{\sum_{\ell=1}^{L}a_\ell}\sum_{\ell=1}^{L}a_\ell
\frac
{m^\ell
}{(m -1)T}\Var \bigl(Z^{m^{\ell},m^{\ell-1}}_{T,1} \bigr).
\nonumber
\end{eqnarray}
Otherwise, since $\mathbb P(X_T\notin\mathcal D_f)=0$, applying the
Taylor expansion theorem twice we get
\begin{eqnarray*}
&& f\bigl( X^{\ell,m^{\ell}}_{T}\bigr)-f\bigl( X^{\ell,m^{\ell-1}}_{T}
\bigr)
\\
&&\qquad  = \nabla f(X_T).U^{m^{\ell},m^{\ell-1}}_T +
\bigl(X^{\ell,m^{\ell}}_{T}-X_T\bigr) \varepsilon
\bigl(X_T,X^{\ell,m^{\ell}}_{T}-X_T\bigr)
\\
&&\quad\qquad{} -
\bigl(X^{\ell,m^{\ell-1}}_{T}-X_T\bigr)\varepsilon
\bigl(X_T,X^{\ell,m^{\ell-1}}_{T}-X_T\bigr).
\end{eqnarray*}
The function $\varepsilon$ is given by the Taylor--Young expansion, so
it satisfies\break
$\varepsilon(X_T,X^{\ell,m^{\ell}}_{T}-X_T){\mathop{\longrightarrow}\limits _{\ell\to\infty}^{{\mathbb
P}}}0$ and
$\varepsilon(X_T,X^{\ell,m^{\ell-1}}_{T}-X_T)
 {\mathop{\longrightarrow}\limits _{\ell\to\infty}^{{\mathbb
P}}}0$. By
property (\ref{Pbold}), we
get the tightness of $\sqrt{\frac{m^{\ell}}{(m-1)T}}(X^{\ell,m^{\ell
}}_{T}-X_T)$ and $\sqrt{\frac{ m^{\ell}}{(m-1)T}}(X^{\ell,m^{\ell
-1}}_{T}-X_T)$
and then we deduce
\begin{eqnarray*}
&& \sqrt{\frac{ m^{\ell}}{(m-1)T}} \bigl(\bigl(X^{\ell,m^{\ell
}}_{T}-X_T
\bigr)\varepsilon\bigl(X_T,X^{\ell,m^{\ell}}_{T}
-X_T\bigr)
\\
&&\hspace*{56pt}{} - \bigl(X^{\ell,m^{\ell-1}}_{T}
-X_T\bigr)\varepsilon\bigl(X_T,X^{\ell,m^{\ell-1}}_{T}-X_T
\bigr) \bigr)  {\mathop{\longrightarrow} _{\ell\to\infty}^{{\mathbb P}}}0.
\end{eqnarray*}
So, according to Lemma \ref{lemma} and Theorem \ref{th-acc}
we conclude that
%e21 #&#
%
\begin{eqnarray}\label{cle-euler}
\sqrt{\frac{m^{\ell}}{(m-1)T}} \bigl(f\bigl( X^{\ell,m^{\ell}}_{T}
\bigr)-f\bigl( X^{\ell,m^{\ell-1}}_{T}\bigr) \bigr) \Rsta
\nabla f(X_T).U_T
\nonumber\\[-10pt]\\[-10pt]
\eqntext{\mbox{as } \ell\rightarrow\infty.}
\end{eqnarray}
Using (\ref{Hf}) it follows from property (\ref{Pbold}) that
\[
\forall\varepsilon>0\qquad \sup_{\ell}\mathbb E\biggl\llvert \sqrt{
\frac
{m^{\ell}}{(m-1)T}} \bigl(f\bigl( X^{\ell,m^{\ell}}_{T}\bigr)-f\bigl(
X^{\ell,m^{\ell
-1}}_{T}\bigr) \bigr) \biggr\rrvert ^{2+\varepsilon}<
\infty.
\]
We deduce using relation (\ref{cle-euler}) that
\begin{eqnarray}
\mathbb E \biggl(\sqrt{\frac{ m^{\ell}}{(m-1)T}} \bigl(f\bigl( X^{\ell,m^{\ell
}}_{T}
\bigr)-f\bigl( X^{\ell,m^{\ell-1}}_{T}\bigr) \bigr) \biggr)^k
\rightarrow \widetilde {\mathbb E} \bigl(\nabla f(X_T).U_T
\bigr)^k<\infty\nonumber
\\
\eqntext{\mbox{for } k\in\{1,2\}.}
\end{eqnarray}
Consequently,
\[
\label{var-euler} \frac{ m^{\ell}}{(m-1)T} \Var\bigl( Z^{m^{\ell},m^{\ell-1}}_{T,1}
\bigr) \rightarrow\widetilde \Var \bigl(\nabla f(X_T).U_T
\bigr)<\infty.
\]
Hence, combining this result together with relation (\ref
{TCL-LindebergFeller}), we obtain the first condition using
Toeplitz lemma. Concerning the second one, by Burkh\"older's inequality
and elementary computations, we get for $p>2$
%e22 #&#
%
\begin{equation}
\label{eq-momentinequality} \mathbb E|X_{n,\ell}|^p=\frac{n^p}{N_\ell^p}
\mathbb E\Biggl\llvert \sum_{\ell=1}^{N_\ell}
Z^{m^{\ell},m^{\ell-1}}_{T,1}\Biggr\rrvert ^p\leq C_p
\frac{n^p}{N_\ell
^{p/2}}\mathbb E\bigl\llvert Z^{m^{\ell},m^{\ell-1}}_{T,1}\bigr
\rrvert ^p,
\end{equation}
where $C_p$ is a numerical constant depending only on $p$.
Otherwise, property (\ref{Pbold})
ensures the existence of a constant $K_p>0$ such that
\[
\mathbb E \bigl| Z^{m^{\ell},m^{\ell-1}}_{T,1} \bigr|^p \leq
\frac{K_p}{m^{p\ell/2}}.
\]
Therefore,
%e23 #&#
%
\begin{eqnarray}\label{eq-momentinequalitybis}
\sum_{\ell=1}^{L}\mathbb E
\llvert X_{n,\ell}\rrvert ^{p} &\leq&\widetilde C_p \sum
_{\ell=1}^{L} \frac{n^p}{N_\ell^{p/2}m^{p\ell/2}}
\nonumber\\[-8pt]\\[-8pt]
&\leq&
\frac{\widetilde C_p}{ (\sum_{\ell=1}^{L}a_\ell )^{p/2}}\sum_{\ell
=1}^{L}a_\ell^{p/2}
\mathop{\longrightarrow}_{n\rightarrow\infty} 0.\nonumber
\end{eqnarray}
This completes the proof.
\end{pf*}

%re1 #&#
%
\begin{remark}
\label{rem}
From Theorem 2, page 544 in \cite{Fel}, we prove a Berry--Esseen-type
bound on our central limit theorem. This improves the relevance of the
above result.
Indeed, take $\alpha=1$ as in the proof, for $X_{n,0}=n\widehat Q_n^1$ and
$X_{n,\ell}$ given by relation~(\ref{Xnl-BerryESSen}), with $\ell\in
\{
1,\ldots,L\}$,
put
\[
s_n^2=\sum_{\ell=0}^{L}
\mathbb{E} |X_{n,\ell}|^{2},\qquad\rho _n=\sum
_{\ell=0}^{L} \mathbb{E} |X_{n,\ell}|^{3}
\]
and denote by $F_n$ the distribution function of $n(Q_n-\mathbb
Ef(X^n_T))/s_n$.
Then for all $x\in\mathbb R$ and $n\in\mathbb N^*$
%e24 #&#
%
\begin{equation}
\label{eq-BerryEssen} \bigl\llvert F_n(x)-G(x)\bigr\rrvert \leq6
\frac{\rho_n}{s_n^3},
\end{equation}
where $G$ is the distribution function of a standard Gaussian random
variable. If we interpret the output of the above
inequality as sum of independent individual path simulation, we get
\begin{eqnarray*}
s_n^2&=&\frac{1}{(m-1)T\sum_{\ell=1}^{L}a_\ell}
\\
&&{}\times \Biggl( a_0 \Var
\bigl(f\bigl( X^1_{T}\bigr) \bigr)+ \sum
_{\ell=1}^{L}a_\ell m^{\ell}\Var
\bigl(f\bigl( X^{\ell,m^{\ell}}_{T}\bigr)-f\bigl( X^{\ell,m^{\ell
-1}}_{T}
\bigr) \bigr) \Biggr).
\end{eqnarray*}
According to the above proof, it is clear that $s_n$ behaves like a
constant but getting lower bounds for $s_n$ seems not to be a common
result to our knowledge.
Concerning $\rho_n$, taking $p=3$ in both inequalities (\ref
{eq-momentinequality}) and (\ref{eq-momentinequalitybis}) gives us an
upper bound. In fact, when $f$ is Lipschitz, there exists a positive
constant $C$ depending on $b$, $\sigma$, $T$ and $f$ such that
\[
\rho_n\leq\frac{ C}{ (\sum_{\ell=1}^{L}a_\ell
)^{3/2}}\sum_{\ell=1}^{L}a_\ell^{3/2}.
\]
For the optimal choice $a_\ell=1$, given in the below subsection, the
obtained Berry--Esseen-type bound is of order $1/\sqrt{\log n}$.
\end{remark}

%
%re2 #&#
%
\begin{remark}
\label{rem-keb}
Note that the above proof differs from the ones in Kebaier
\cite{Keb}. In fact, here our proof is based on the central limit
theorem for triangular array which is adapted to the form of the
multilevel estimator, whereas Kebaier used another approach based on
studying the associated characteristic function. Further, this latter
approach needs a control on the third moment, whereas we only
need to control a moment strictly greater than two. Also, it is worth
to note that the limit variance in Theorem \ref{CLTSR} is smaller than the
limit variance in Theorem 3.2 obtained by Kebaier in \cite{Keb}.
\end{remark}

%s3.3 #&#
\subsection{Complexity analysis}\label{euler-parameters}

\begin{figure}%[hbt]

\includegraphics{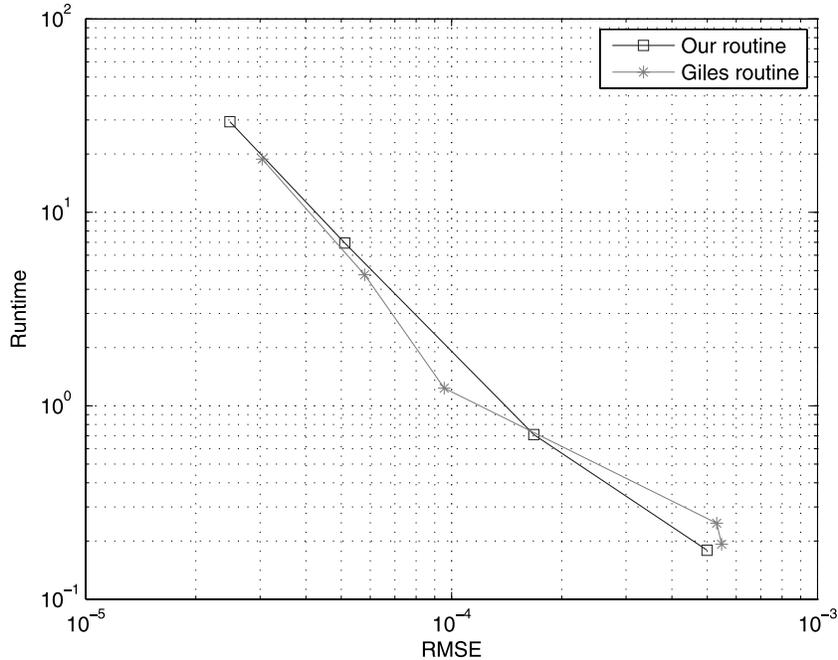}

\caption{Comparison of both routines.}\label{fig1}
\end{figure}

From a complexity analysis point of view, we can interpret Theorem \ref
{CLTSR} as follows.
For a total error of order $1/n^{\alpha}$, the computational effort
necessary to run the multilevel Monte Carlo Euler method
is given by the sequence of sample sizes specified by relation (\ref{eq-size}).
The associated time complexity is given by
\begin{eqnarray*}
C_{\MMC} &= &C\times \Biggl(N_0 + \sum
_{\ell=1}^{L}N_{\ell
}\bigl(m^{\ell}+m^{\ell-1}
\bigr) \Biggr)\qquad\mbox{with } C>0
\\
&=& C\times \Biggl( \frac{n^{2\alpha}(m-1)T}{a_0}\sum_{\ell=1}^{L}a_\ell+
n^{2\alpha}\frac{(m^2-1)T}{m}\sum_{\ell=1}^{L}
\frac{1}{a_\ell
}\sum_{\ell
=1}^{L}a_\ell
\Biggr).
\end{eqnarray*}
The minimum of the second term of this complexity is reached for the
choice of weights $a_\ell^*= 1$, $\ell\in\{1,\ldots,L\}$,
since the Cauchy--Schwarz inequality ensures that $L^2\leq\sum_{\ell
=1}^{L}\frac{1}{a_\ell}\sum_{\ell=1}^{L}a_\ell$, and the
optimal complexity for the multilevel\vspace*{2pt}
Monte Carlo Euler method is given by
\begin{eqnarray*}
\label{com-euler}
C_{\MMC} &=& C\times \biggl( \frac{(m-1)T}{a_0 \log m}
n^{2\alpha} \log n+\frac{(m^2-1)T}{m(\log m)^2 } n^{2\alpha}(\log n)^2
\biggr)
\\
&=& \mathit O \bigl(n^{2\alpha}(\log n)^2 \bigr).
\end{eqnarray*}
It turns out that for a given discretization error $\varepsilon
_n=1/n^\alpha$ to be achieved the complexity is given by
$C_{\MMC}=\mathit O ( \varepsilon_n^{-2}(\log\varepsilon
_n)^2 )$.
Note that this optimal choice $a_\ell^*= 1$, $\ell\in\{1,\ldots,L\}$,
with taking $a_0=1$ corresponds to the sample sizes given by
\[
N_\ell=\frac{(m-1)T }{m^{\ell}\log m}n^{2\alpha}\log n,\qquad \ell\in \{ 0,\ldots,L\}.
\]
Hence, our optimal choice is consistent with that proposed by Giles
\cite{Gil}.
Nevertheless, unlike the parameters obtained by Giles \cite{Gil} for
the same setting [see relation~(\ref{GilesParam})], our optimal choice of
the sample sizes $N_{\ell}, \ell\in\{1,\ldots,L\}$ does not depend on
any given constant, since our approach is based on proving a central
limit theorem and not on getting upper bounds for the variance.
Otherwise, for the same error of order $\varepsilon_n=1/n^{\alpha}$ the
optimal complexity of a Monte Carlo method is given by
\[
C_{\MC}= \mathit O \bigl( n^{2\alpha+1} \bigr)= \mathit O \bigl(
\varepsilon_n^{-2 -1/\alpha} \bigr)
\]
which is clearly larger than $C_{\MMC}$. So, we deduce that the multilevel
method is more efficient. Also, note that the optimal choice of
the parameter $m$ is obtained for $m^*=7$. Otherwise, any choice
$N_0=n^{2\alpha}(\log n)^{\beta}$, $0<\beta<2$, leads to the same result.
Some numerical tests comparing original Giles work \cite{Gil} with the
one of us show that both error rates are in line. Here in Figure~\ref{fig1}, we
make a simple log--log scale plot of CPU time with respect to the root
mean square error, for
European call and with $N_0=n^{2\alpha}(\log n)^{1.9}$.\vadjust{\goodbreak}

%
%f1 #&#
%

It is worth to note that the advantage of the central limit theorem is
to construct a more accurate confidence interval. In fact, for a given
root mean square error RMSE, the radius of the $90\%$-confidence
interval constructed by
the central limit theorem is
$1.64\times{}$RMSE. However, without this latter result, one can only use
Chebyshev's inequality which yields a radius equal to $3.16\times{}$RMSE.
Finally, note that, taking $\alpha=1/2$ still gives the optimal rate
and allows us to cancel the bias in the central limit theorem due to
the Euler discretization.

%s4 #&#
\section{Conclusion}
The multilevel Monte Carlo algorithm is a method that can be used in a
general framework:
as soon as we use a discretization scheme in order to compute
quantities such as $\mathbb Ef ( X_T )$,
we can implement the statistical multilevel algorithm. And this is
worth because it is an efficient method according to the original work
by Giles \cite{Gil}.
The central limit theorems derived in this paper fill the gap in
literature and confirm superiority of the multilevel method over the
classical Monte Carlo approach.

%pa4.subsection.subsubsection.1 #&#
% zodis "Acknowledgments" paliekamas pagal autoriu
\section*{Acknowledgments}
The authors are greatly indebted to the Associate Editor and the
referees for their many suggestions and interesting comments which
improve the paper.

%suskaldyti doi

% imsref loaded by linak, 2014-03-21 11:21:13
%

\printaddresses

\end{document}